\documentclass[11pt]{article}

\usepackage[utf8]{inputenc}
\usepackage[T1]{fontenc}
\usepackage{geometry}
\usepackage{amsmath}
\usepackage{amssymb}
\usepackage{amsthm}
\usepackage{mathtools}
\usepackage{graphicx}
\usepackage{mathrsfs}
\usepackage{xcolor}
\usepackage{booktabs}
\usepackage{caption}
\usepackage{enumitem}
\usepackage{braket}
\usepackage[round]{natbib}
\usepackage{comment}
\usepackage[hidelinks]{hyperref}
\usepackage{doi}
\usepackage{microtype}
\usepackage{cleveref}

\ifdefined\draftmode
\usepackage{draftwatermark}
\SetWatermarkText{DRAFT}
\SetWatermarkScale{3}
\newcommand{\tk}[1]{\textcolor{red}{\textbf{[Neil: #1]}}}
\else
\newcommand{\tk}[1]{}
\fi

\theoremstyle{definition}
\newtheorem{definition}{Definition}
\theoremstyle{remark}

\numberwithin{theorem}{section}
\numberwithin{proposition}{section}
\numberwithin{lemma}{section}
\numberwithin{corollary}{section}
\numberwithin{definition}{section}
\numberwithin{remark}{section}

\ifdefined\draftmode
\else
\excludecomment{draftcomment}
\fi

\begin{document}

\title{The No Barber Principle: Towards Formalised Selection in the Inaccessible Game}

\author{Neil D. Lawrence\\
Computer Laboratory\\
University of Cambridge\\
15 J. J. Thomson Avenue\\
Cambridge CB3 0FD, UK\\
\href{mailto:ndl21@cam.ac.uk}{ndl21@cam.ac.uk}}

\date{21st April 2026}

\maketitle

\begin{abstract}
  The inaccessible game
  \citep{Lawrence-inaccessible25,Lawrence-origin26} is an
  information-theoretic dynamical system governed by three information
  loss axioms, a marginal entropy conservation constraint and maximum
  entropy dynamics. In this paper we look at \emph{selection} in the
  game. Our aim is to develop a selection policy for the game rules
  based on a minimal set of assumptions. We seek necessary consistency
  constraints for self-determining dynamical systems. Specifically, we
  suggest that rules that quantify over distinctions they cannot
  internally represent risk impredicative-style circularity. Our
  criterion is motivated by an analogy with Russell's paradox. We
  formulate a \emph{no-barber principle} which prohibits dynamics that
  appeal to external adjudicators or structure lying outside the
  system.

  To motivate our principle we examine Russell's paradox through its
  structural formalisation as a Lawvere diagonalisation. The
  marginal-entropy conservation in the game is a nontrivial entropy
  constraint which prohibits external structure.   Through the no-barber
  principle we argue (i) the classical category \textsf{FinProb}, in
  which Shannon entropy is characterised, is cartesian and provides
  canonical diagonal (copying) maps that make Lawvere-style
  constructions expressible and is structurally incompatible with the
  no-copying instantiation of the no-barber principle studied
  here. (ii) the noncommutative category \textsf{NCFinProb}, in which
  von Neumann entropy is characterised, is symmetric monoidal and
  lacks canonical copying maps, making it a more natural candidate for
  the game's internal language.
\end{abstract}

\section{Introduction}

The inaccessible game \citep{Lawrence-inaccessible25} is an
information-theoretic game inspired by Conway's Game of Life and
Wolfram's one-dimensional cellular automata. In both those games the
rules of the automata are explicitly defined, but for the inaccessible
game we would like the rules to be emergent.

\citet{Lawrence-inaccessible25,Lawrence-origin26} explored the
information geometry of the game and derived some promising
mathematical structure. However, to do this explicit design
choices were made. For example, in \citet{Lawrence-inaccessible25} the
game is assumed to proceed through maximum entropy dynamics. In
\citet{Lawrence-origin26} the origin point is assumed to be located
where the marginal entropies are maximised. When justifying these
design choices \citet{Lawrence-origin26} introduced an informal notion
of `axiomatically distinguished'. The idea being to ensure that any
design selections were specified uniquely within the internal
mathematical structure of the game.

This paper is the first step in formalising how such axiomatically
distinguished choices might be made. To motivate our approach consider
the card game Munchkin \citep{Jackson-munchkin01}. It has written
rules to deal with inconsistency in game rules. The suggested
resolution is that `... disputes should be settled by loud arguments,
with the owner of the game having the last word.' Such external
adjudication provisions are common across sports, board games and even
societal rules. But this raises the question: is this Munchkin
provision always necessary? Or can we constrain a game to operate
without external adjudication?

We consider a proposal that aims to eliminate conditions where
external adjudication is necessary. Our proposal doesn't justify all
the choices made in the game so far, but it does give insight into why
the fundamental category underpinning the information loss cannot be
cartesian, eliminating the axiomatic form of information loss that
originally motivated the game.

Let's first consider an example of a system where external
adjudication would be required. Russell's paradox
\citep{Russell-frege02,Deutsch-russell24} is foundational in the idea
of impredicative circularity. The canonical description of this
paradox involves a barber who shaves everyone in the village that
doesn't shave themselves. The impredicative circularity is associated
with the question: `Does the barber shave themself?'. If such a
question needed resolving for the dynamics of a game to proceed then
it would require external adjudication. We therefore need a selection
principle that prevents the possibility of impredicative
circularity. In honour of the canonical example of Russell's paradox,
we call this principle the no-barber principle because it aims to
remove the impredicative circularity associated with the barber.
\begin{definition}[No-Barber Principle (informal)]
An internal language is admissible for a self-adjudicating dynamical
system only if it does not permit constructions whose resolution
requires appeal to external adjudication or structure that lies
outside the system's own dynamics.
\end{definition}
We study one sufficient structural proxy for this principle: the
absence of canonical copying maps.

First we note that the information isolation axiom
\citep{Lawrence-inaccessible25} is the foundation of the game. This
should prevent any access to external structure or adjudication. It
implies that the game is built `without introducing any additional
external choice such as Hamiltonians, spatial coordinates, a
temperature scale, or a background clock.' But beyond the blocking
implied by the information isolation axiom, current design choices are
informal. The information isolation axiom is based on axioms that
define information loss through category theory
\citep{Baez-characterisation11,Parzygnat-functorial22}. This raises a
challenge, how do we know that our definition of information loss is
not smuggling in some form of external structure? We will examine this
question using the no-barber principle.

To better formalise the principle, we rely on William Lawvere's
characterisation of Russell's paradox and a broader class of related
paradoxes through a diagonalisation argument
\citep{Lawvere-diagonal69}. This allows us to identify that the
structure of diagonalisation can be blocked through preventing
cartesian copying. We then suggest that the elimination of canonical
copying maps is a \emph{design constraint} that is
consistent with the no-barber principle. It eliminates a class of impredicative
circularity for self-adjudicating dynamical systems.\footnote{Note that the restriction is not a general restriction on
mathematical formalisms we use for analysis or even on emergent
structures within the game. It is a restriction on the game's internal
dynamics governing the system's own evolution.}

%% \citet{Lawvere-diagonal69}'s diagonalisation argument applies to a
%% wider class of methods where impredicative circularity can
%% emerge. Alongside Russell's paradox the same structure underpins
%% G\"odel's incompleteness theorem \citep{Godel-incompleteness31},
%% Turing's proof for the halting problem \citep{Turing-computable36} and
%% Tarski's undefinability theorem \citep{Tarski-warhheitsbegriff36}.

%% This unifying perspective suggests that we should seek to prohibit
%% Lawvere's structure from forming in the internal dynamics. Blocking
%% the formation gives us a selection principle that eliminates these
%% classes of impredicative circularity from the game rules.

The rest of the paper is structured as
follows. Section~\ref{sec:diagonal-mechanism} reviews Lawvere
diagonalisation and formulates the no-barber principle as a design
constraint that blocks copying and hence the diagonal
construction. Then in Section~\ref{sec:linear-logic} we highlight that
existing powerful frameworks including linear logic and quantum
probability instantiate the relevant design constraint. We conclude in
Section~\ref{sec:von-neumann} by noting that this design constraint is
consistent with \textsf{NCFinProb} (and von Neumann entropy) over
\textsf{FinProb} (and Shannon entropy) as the internal language for
the game, a choice already independently motivated in the origin paper
\citep{Lawrence-origin26}. The paper closes with a discussion and
future-work section (Section~\ref{sec-discussion}).

\section{The Diagonalisation Mechanism}
\label{sec:diagonal-mechanism}

In this section we introduce the structure of Lawvere's
diagonalisation. The arguments come from category theory, although
here we will motivate only by the structure of the argument leaving a
full rigorous category theoretic formulation to further work.

From a categorical perspective \citep{Lawvere-conceptual09} we
consider the rules of our game to be elements of an internal language
whose evaluation is represented by morphisms in a chosen category. The
structural properties of that category then determine which
constructions, such as duplication and self-application, are
available.

We've motivated this paper by arguing that for the inaccessible game
\citep{Lawrence-inaccessible25,Lawrence-origin26} these rules need to
be self-adjudicating. How can we constrain the game so that it has
this property?

By inheriting a notion of information loss from
\citet{Baez-characterisation11} and \citet{Parzygnat-functorial22}),
the game expresses a notion of information conservation. But how can
we be sure that the underlying definition of information loss does not
smuggle in external structure? In the introduction we suggested that
avoiding impredicative circularity may provide the key. Motivated by
this idea we briefly review Lawvere diagonalisation
\citep{Lawvere-diagonal69,Lawvere-conceptual09}, which provides a
structure through which self-referential systems fail in impredicative
circularity. We will then show one route to preventing this structure
from forming by prohibiting cloning within our system. To illustrate
the idea we further review Russell's paradox through the lens of
Lawvere's diagonalisation.

In 1901, when formulating the foundations of mathematics, Bertrand
Russell identified a paradox via a set $R$ that is defined as the set
of all sets that are not members of themselves
\citep{Deutsch-russell24}. This triggers the question: is $R \in R$?
If so then by definition $R \notin R$. Conversely if $R \notin R$ then
$R \in R$. Russell illustrated this paradox by defining a set of
villagers as shaved by a barber only if they don't shave
themselves. Because the barber is a member of the village the
structure of the paradox is triggered.

The paradox is the canonical example of \emph{impredicative
definition}. It arises because we are quantifying over a totality that
includes the object being defined. If we translate this to the game
dynamics, it's equivalent to a rule that says `do $X$ based on
checking a property of all objects,'. But then that property is in
turn dependent on applying the rule. Under these conditions
impredicative circularity can arise. Such an unresolvable rule in the
dynamics of the game would require external arbitration to prevent
stalling in circular reference.

How can we block such structure in the game?
\citet{Lawvere-diagonal69} showed that Russell's barber paradox has a
common category theoretic structure with other proofs including
Cantor's theorem \cite{Cantor-elementare91}, G\"odel incompleteness
\citep{Godel-incompleteness31}, Turing's halting problem
\citep{Turing-computable36}, and Tarski's undefinability theorem
\citep{Tarski-warhheitsbegriff36}. We informally summarise this
structure as a five-move recipe.
\begin{enumerate}[label=\arabic*.]
\item The system consists of \emph{rules} which are a collection of
  operations that take inputs to outputs. \label{lawvere-rules}
\item There must be a \emph{universal evaluator} operation that can
  execute a rule on any input. \label{lawvere-evaluator}
\item For every rule it must be possible to \emph{internalise the
syntax}, i.e.\ the system should be able to identify all its operations
  internally. \label{lawvere-internalise}
\item The rules should allow for \emph{self-application}, so that any
  rule-code can be fed in as an input, i.e.\ the system's rules can
  serve as their own input. \label{lawvere-self}
\item There should be a \emph{twist} transformation that can negate or
  complement outputs. \label{lawvere-twist}
\end{enumerate}
This structure applies across an impressive range of mathematical
results. They illustrate that we can't tell whether a (Turing complete)
computer programme will stall. Or whether a given mathematical
statement that is true is provable (G\"odel). Or even whether truth itself is internally definable \citep{Tarski-warhheitsbegriff36}.

To apply the no-barber principle, we are looking to break this chain
at some point. The overall structure provides several different points
for intervention. For example self-application
(Step~\ref{lawvere-self}) requires feeding a rule-code in as its own
input. This corresponds internally to duplicating that input: the same
object must appear simultaneously as both the rule and its
argument. Lawvere develops his proof in a cartesian setting where this
operation is represented by the canonical diagonal map $A \to A \times
A$. However, in a monoidal setting there is no corresponding canonical
copying map $A \to A \otimes A$.\footnote{In a cartesian category the
product $\times$ is equipped with canonical diagonal maps
$\Delta_A: A \to A \times A$ as part of the categorical
structure. The tensor product $\otimes$ in a symmetric monoidal
category carries no such universal diagonal. Any copying must be
explicitly postulated.} In the absence of such a map,
Step~\ref{lawvere-self} cannot be formed and the construction fails.
Note that we are not claiming that every violation of the no-barber
principle arises from copying, but noting that canonical copying is a
structural source of the relevant impredicative failure, and hence a
natural target for exclusion.

\begin{quote}
\textbf{Proposition (informal).}  In any internal language lacking a
canonical copying map $A \to A \otimes A$, the self-application step
required for Lawvere diagonalisation cannot be formed without
additional, non-standard structure.  Consequently, the standard
Lawvere-style fixed-point contradiction, including the barber paradox
and related diagonal arguments, cannot be derived internally.
\end{quote}

When all five moves are available, one can construct a diagonal rule
of the form $g(x) = \text{Twist}(\text{Apply}(x, x))$. Through
internalisation of syntax, there must exist some internal code $r_0$
such that we can form $g(x) = \text{Apply}(r_0, x)$ for all
$x$. Evaluating at $x = r_0$ gives the fixed point: $\text{Apply}(r_0,
r_0) = \text{Twist}(\text{Apply}(r_0, r_0))$, which is a contradiction
for appropriate choices of twist like negation or complement.

Being categorical in character, this argument is not specific to sets,
proofs, or programs. Any mathematical system admitting all five moves
can exhibit paradox or incompleteness. As an example let's look at how
the structure applies to the case of the barber paradox.

Let $P$ be the set of people in the village and $R$ the set of shaving
rules (stage \ref{lawvere-rules}),
$$
R : P \to \{\text{yes},\text{no}\}
$$
that specify which people should be shaved.

Combining a rule with a person produces a shaver object, $S = R \times
P$
$$
x=(R,p) \in S.
$$
This represents the person $p$ operating according to rule $R$.

There is a universal evaluator operation (stage \ref{lawvere-evaluator})
$$
\text{shaves} : S \times P \to \{\text{yes},\text{no}\}
$$
which applies the rule stored in a shaver object to a person
$$
\text{shaves}((R,q),p) = R(p).
$$
Thus any rule embedded in a shaver object can be executed on any
person.

The internalisation of the rules (stage \ref{lawvere-internalise})
arises because each person $p$ has some shaving behaviour. We
associate a rule, $R_p$, describing that behaviour. So each person
determines a canonical shaver object consisting of a person and their
shaving rule
$$
x_p=(R_p,p).
$$
This means that when a rule receives a person $p$ as input, it can
recover and evaluate the rule governing that person's behaviour. That
implies that the rule must be internalised.

Since rules act on people, we can apply a person's rule to themself,
$R_p(p)$ (stage \ref{lawvere-self}). In terms of the evaluator this is
$$
\text{shaves}(x_p,p).
$$ Using this capability we now define the twist rule (stage
\ref{lawvere-twist}),
$$
T(p)=\neg R_p(p).
$$
This rule takes a person $p$, retrieves their rule $R_p$, applies it
to $p$, and negates the result. Thus $T$ selects exactly those people
who do not shave themselves.

We now suppose there were a person $b$ whose shaving behaviour
followed this rule. This gives us the barber who would be,
$$
B=(T,b).
$$
Evaluating the barber on themself gives
$$
\text{shaves}(B,b)=T(b)=\neg R_b(b).
$$
But since $R_b=T$, this becomes
$$
\text{shaves}(B,b)=\neg \text{shaves}(B,b).
$$
This gives us the contradiction in the paradox. The relationship
between each of these steps and Lawvere's rules is given in Table
\ref{table-barber-lawvere}.

\begin{table}
\caption{Correspondence between Lawvere's steps and the barber
  paradox.\label{table-barber-lawvere}}
\begin{center}
  \begin{tabular}{ll}
    \hline
    Step & Barber interpretation \\
    \hline
    rules & shaving rules $R$ \\
    universal evaluator & $\text{shaves}((R,p),q)$ \\
    internalisation	& each person $p$ carries rule $R_p$ \\
    self-application & $R_p(p)$ \\
    twist & $T(p)=\neg R_p(p)$ \\
    \hline
  \end{tabular}
  \end{center}
\end{table}
The barber paradox therefore arises because the system allows rules to
inspect and invert their own behaviour via the self-application
$R_p(p)$.

\subsection{Blocking the Diagonalisation}

The no-barber principle suggests that we need to block the formation
of these Lawvere diagonalisations to prevent impredicative
circularity. There are various possible points of intervention, but we
will look specifically at the diagonalisation itself.

Diagonalisation in the categorical sense means the move of restricting
a \emph{two}-input evaluation map to the diagonal by duplicating its
input. Concretely, in a category with finite products, every object
$P$ has a canonical diagonal morphism $\Delta_P \colon P \to P \times
P$. Given a `bivariate' morphism $F \colon P \times P \to Y$, or
more generally an evaluator $\text{Apply}\colon A\times P\to Y$
together with a way of viewing an element of $P$ as determining a code
in $A$, the diagonalised map is the composite $ P
\xrightarrow{\Delta_P} P\times P \xrightarrow{F} Y, $ i.e.\ $p \mapsto
F(p,p)$. This is a `feed the same person twice' operation: one copy of
$p$ is used in the `rule' slot and the other in the `input'
slot. Lawvere-style paradoxes then arise when we post-compose this
diagonal evaluation with a twist $\alpha \colon Y\to Y$, such as
negation in the Barber case, to form $p \mapsto \alpha(F(p,p))$.

This description highlights the duplication assumption in the
diagonalisation. To compute the self-application $R_p(p)$ appearing in
the twist rule $T(p)=\neg R_p(p)$, the system must use the same person
$p$ in two roles simultaneously: one copy is used to recover the rule
$R_p$ governing that person's shaving behaviour, while the other copy
is supplied as the argument to that rule. The duplication step is
implicit in the barber construction, but is essential for forming the
diagonal rule. Without this cloning operation the twist rule $T(p)$
could not be defined, and the diagonal contradiction would not arise.

This suggests that if we prevent cloning in the system we block the
diagonalisation. We could also try to block at another point in the
chain, but among the five moves we view blocking copying as the most
structurally natural intervention for three reasons. First, in the
diagonalisation copying is not introduced as an explicit rule but is
enforced for free by cartesian structure. Second, removing it does not
ban particular rules such as the twist or the universal evaluator. It
alters the ambient structure, leaving expressive power otherwise
intact.  Third, as we'll review in the next section, the formalisms
of linear logic and quantum already have prohibitions of canonical
copying maps demonstrating the existence of well-developed,
expressive formal systems where structural copying is already prohibited.

\subsection{Linear Logic and Quantum Probability}
\label{sec:linear-logic}

The suggestion `block diagonalisation by blocking copying' has a
useful interpretive parallel with two reconstructionist programmes.
We now briefly review two frameworks that prohibit copying.

Ed Jaynes \citep{Jaynes-probability03} argues that ordinary
probability theory is best viewed as an extension of classical logic
to plausible reasoning under incomplete information. Once we assume an
underlying Boolean event structure (with conjunction, disjunction,
negation) and demand consistency of inference, then we are effectively
forced into the usual probability calculus. In that classical setting,
the structural rules of reasoning are invisible but strong: premises
can be reused arbitrarily and discarded freely. There is no block on
copying.

In converse, Girard's linear logic \citep{Girard-linear87} makes those hidden
structural rules explicit by removing them. In linear logic,
assumptions behave like resources that cannot in general be duplicated
or discarded. Categorically this corresponds to moving from a
cartesian world, where every object has a canonical diagonal $A\to
A\times A$ and terminal map $A\to 1$, to a symmetric monoidal world
$(\otimes,I)$, where there is in general no canonical copy map $A\to
A\otimes A$. Copying is reintroduced only for distinguished
`of-course' objects, where a chosen comonoid structure provides
controlled duplication.

Quantum probability fits naturally into this resource-sensitive
picture. The algebraic approach replaces Boolean event algebras by
noncommutative $C^*$-algebras of observables, and replaces functions
between sample spaces by structure-preserving maps between algebras
\citep{Petz-quantum08,Parzygnat-functorial22}. In this setting there
is no state-independent way to duplicate an arbitrary unknown
state. This is typically referred to as `no-cloning'. In the classical
(commutative) case diagonals are canonical because they come from
copying points of a sample space. In the noncommutative case that
copying map fails.\footnote{This begs the question, beyond the scope
of this work, as to whether the shift from classical to noncommutative
probability is analogous (in Jaynes's programme) to replacing ordinary
Boolean logic by a compositional, resource-sensitive logic.}

The no-barber principle suggests that Lawvere-style diagonal
constructions are not an unavoidable feature of expressive systems,
but a design choice that we reject. This imposes a design constraint
on our internal language. Linear logic and quantum probability are two
established and expressive programmes that block copying. They are
naturally expressed in symmetric monoidal terms ensuring that the
diagonal map is not available by default.

\subsection{Von Neumann Entropy}
\label{sec:von-neumann}

The inaccessible game relies on characterising information loss so
that it can impose information isolation. In
\citet{Parzygnat-functorial22} information loss is characterised
through von Neumann entropy in the category of finite-dimensional
noncommutative probability spaces (\textsf{NCFinProb}). In
\citet{Baez-characterisation11} information loss is characterised
through Shannon entropy in the classical category \textsf{FinProb}.

In \textsf{FinProb}, copying is canonical: the cartesian diagonal maps
$A \to A \times A$ are available for free, and with them the diagonal
constructions that underwrite impredicative paradox. In
\textsf{NCFinProb}, a symmetric monoidal category, canonical copying
does not exist. Under our no-copying proxy the no-barber principle
counts against \textsf{FinProb} as the foundational internal
language for a self-adjudicating system and points towards
\textsf{NCFinProb} as a more natural candidate for the game's internal
language.

This choice is independently motivated in \citet{Lawrence-origin26},
where it is shown that the origin of the game requires zero joint
entropy with positive marginal entropies. This is a configuration that
is impossible under Shannon entropy (whose conditional entropy is
always non-negative) but achievable under von Neumann entropy through
entanglement. The no-barber principle reinforces the selection made in
that paper but does so through independent structural grounds.

\section{Discussion and Future Work}
\label{sec-discussion}

The objective of this paper was to provide a selection principle for
the inaccessible game. We have argued that any self-adjudicating
dynamical system must have an internal language that prohibits
stalling in impredicative circularity. We called this the no-barber
principle. Through Lawvere's analysis of impredicative circularity we
suggested a structural conclusion: the internal language of such a
system should not admit canonical copying maps. As a result it should
not be grounded in a cartesian category.

For the inaccessible game this provides a structural reason not to
take \textsf{FinProb} as the foundational setting for information
loss, and points to \textsf{NCFinProb}, and the von Neumann
entropy, as the more natural alternative.

Our arguments are informal, but suggestive. What a fully rigorous
treatment would require is formal definitions of admissible internal
languages, a precise categorical characterisation of canonical
copying. These challenges are left to future work.

We emphasise that we are \emph{not} claiming to evade G\"odelian
limitations or propose a new foundation for mathematics. The
no-copying restriction applies to the internal language governing the
evolution of the dynamical system. It does not constrain the
meta-mathematical reasoning used to analyse the system. It also does
not prevent subsystems with different effective behaviour from
emerging within the system.

\bibliographystyle{plainnat}
\bibliography{the-inaccessible-game-no-barber}

\end{document}